\documentclass[12pt,twoside]{article}
\usepackage{graphicx}
\usepackage{graphics}
\usepackage{amsmath,amssymb,amsthm}

\theoremstyle{plain}
\newtheorem{thm}{Theorem}[section]

\newtheorem{lem}{Lemma}[section]

\theoremstyle{remark}
\newtheorem{rem}{Remark}[section]

\newtheorem*{pf}{Proof}

\def\E{{\mathbb {E}}}                       \def\D{{\mathbb {D}}}

\def\FD{{\mathcal F}}

\def\o{{\rm {o}}} 
\def\O{{\rm {O}}} 

\begin{document}
\centerline{\large{\bf {Weights of Cliques in a Random Graph Model}}}
\smallskip
\centerline{\large{\bf {Based on Three-Interactions}}}

\bigskip
\centerline{ {\sc Istv\'an Fazekas},  {\sc Csaba Nosz\'aly} and {\sc Attila Perecs\'enyi} }

\medskip
Faculty of Informatics, University of Debrecen,  P.O. Box 12, 4010 Debrecen, Hungary, 
\centerline{
e-mail: fazekas.istvan@inf.unideb.hu.}

\bigskip
\begin{abstract}
A random graph evolution rule is considered.
The graph evolution is based on interactions of three vertices.
The weight of a clique is the number of its interactions.
The asymptotic behaviour of the weights is described.
It is known that the weight distribution of the vertices is asymptotically a power law.
Here it is proved that the weight distributions both of the edges and the triangles are also asymptotically power laws.
The proofs are based on discrete time martingale methods.
Some numerical results are also presented.
\end{abstract}
\renewcommand{\thefootnote}{}
{\footnotetext{


The publication was supported by the T\'AMOP-4.2.2.C-11/1/KONV-2012-0001 project. The project has been supported by the European Union, co-financed by the European Social Fund.
}
%
\section{Introduction}  \label{Sect1}
It is known that many real-life networks are scale free.
Such networks are the WWW, the Internet, certain social and biological networks (see e.g. \cite{durrett}).
A network is called scale-free if it has power law degree distribution, that is $p_k \sim Ck^{-\gamma}$ as $k\to\infty$, 
where $p_k$ is the relative frequency of vertices with degree $k$.
In \cite{barabasi} Barab\'asi and Albert introduced  the preferential attachment rule.
In their model the growing procedure is the following. 
At every time $n=2,3,\dots$ a new vertex is added. 
The new vertex is attached to $m$ existing vertices. 
The probability $\pi_i$ that the new vertex is connected to the old vertex $i$ depends on the degree $k_i$ of that vertex, so that $\pi_i= k_i/\sum_j k_j$.
A rigorous definition of the preferential attachment rule and a precise mathematical proof for the power law degree distribution were given in \cite{bollobas}.

In \cite{BaMo1} a 3-interactions model was introduced.
The evolution in that model is a version of the preferential attachment rule.
The power law degree distribution was obtained in \cite{BaMo2}.
In this paper we shall study the 3-interactions model so we present its detailed definition. 

Our main results will be confined to the following model.
The starting point at time $n=0$ is one triangle.
The initial weight of this triangle is one. 
This graph contains $3$ vertices and $3$ edges. 
Each of them has initial weight $1$.
After the initial step we increase the size of the graph.
At each time step $n=1,2,\dots$ three vertices interact.
The interaction of three vertices means the following.
First we draw all non-existing edges between these vertices.
So we obtain a triangle.
Then the weight of this triangle and the weights of its vertices and edges are increased by $1$. 
Here we should remark that any initial weight is $1$.
That is the weight of a new object (a new triangle, edge or vertex) is $1$.
The choice of the three vertices is the following.

There are two ways at each step.
On the one hand, with probability $p$, we add a new vertex.
I this case the new vertex and two old vertices interact.
On the other hand, with probability $\left( 1-p \right)$, we do not add a new vertex.
In this case three old vertices interact. 
Here $0 < p \leq 1$ is a fixed number.

If we add a new vertex, then we should choose two old vertices and these three vertices will interact.
To choose the two old vertices we have two options.
On the one hand, with probability $r$, we choose an edge from the existing edges according to the weights of the edges.
That is an edge of weight $v_t$ is chosen with probability $v_t/\sum_h v_h$ (it is the preferential attachment rule).
In this case the two end points of the edge chosen play the role of the two old vertices.
On the other hand, with probability $1-r$, we choose two out of the existing vertices uniformly, that is all two vertices have the same chance.

If we do not add any new vertex, then three old vertices interact.
To choose the three old vertices we have again two options.
On the one hand, with probability $q$, we choose from the existing triangles according to their weights.
It means that a triangle with weight $w_t$ is chosen with probability $w_t/\sum_h w_h$ (preferential attachment).
On the other hand, with probability $1-q$, we choose three out of the old vertices uniformly.

For the sake of brevity, a complete graph with $m$ vertices we shall call an $m$-clique.
In the three-interaction model we study the following cliques: vertices, edges and triangles.
Our aim is to show scale free property for the weights of cliques.
For vertices it was proved in \cite{BaMo1}. 
For triangles it was announced in \cite{FaNo}.
Our most interesting result is the power law weight distribution for the edges.
The weight of an edge can be considered as the strength of the connection between its two end points.

The main results are included in Section~\ref{Sect2}.
Theorem~\ref{limEnv/En} is the power law weight distribution of the edges. 
Its proof is included in Section~\ref{Sect3}.
Theorem~\ref{limKnv/Kn} shows the power law distribution of the weights of triangles.
As the result is true for the largest cliques in the general $N$-interactions model, we present the proof in the general case.
Therefore Section~\ref{Sect4} is devoted to the $N$-interactions model.
In the proofs the basic method is the Doob-Meyer decomposition of submartingales.
In Section~\ref{Sect5} computer simulation results are presented.
%
%
\section{Power law distributions of the weights}  \label{Sect2}
\setcounter{equation}{0}
Throughout the paper $0<p \leq 1$, $0 \leq r \leq 1$, $0 \leq q \leq 1$ are fixed numbers.
Let
\begin{equation}  \label{alfabeta}
\alpha= \dfrac{2}{3}pr + \left(1-p\right) q, \qquad  \beta =  2\left( 1-r \right) + \dfrac{3\left( 1-p \right)\left( 1-q \right)}{p}.
\end{equation}
Let $V_n$ denote the number of vertices after the $n$th step.
Let $X \left( n,w \right)$ denote the number of vertices of weight $w$ after $n$ steps.

First we quote Theorem 3.1 of \cite{BaMo1}.
It shows the scale free property for the weights of vertices. 
\begin{thm}  \label{theorem:scalefreeWeights}
Let $0<p<1$, $q>0$, $r>0$.
Assume that either $r<1$ or $q<1$.
Then for all $w=1,2,\dots$ we have
\begin{equation}
\dfrac{X \left( n,w \right)}{V_n} \rightarrow x_{w}
\end{equation}
almost surely as $n \rightarrow \infty $, where $x_{w}$, $w=1,2, \dots$\,, are positive numbers satisfying the recurrence relation
\begin{equation} \label{rekurziox(w)-re}
x_{1} = \dfrac{1}{\alpha + \beta +1},  \qquad
x_{w} = \dfrac{\alpha \left( w-1 \right) + \beta}{\alpha w + \beta +1}x_{w-1}, \quad \text{if} \quad w > 1,
\end{equation}
where $\alpha$ and $\beta$ are defined by \eqref{alfabeta}.
Moreover,
\begin{equation} \label{x_{w}_asz.}
x_{w} \sim C  w^{- \left( 1 + \frac{1}{\alpha} \right)}
\end{equation}
as $w \to \infty$,  with $C =\varGamma \left( 1 + \frac{\beta +1}{\alpha} \right) \big/ \left({\alpha \varGamma \left( 1 + \frac{\beta}{\alpha} \right)}\right) $.
\end{thm}
Here and in what follows $x_n \sim y_n$ means that $\lim_{n\to\infty} x_n/y_n = 1$.
The power law degree distribution is also true (see \cite{BaMo2}).

Let $K(n,w)$ denote the number of triangles with weight $w$, and let $K_n$ denote the number of all triangles after $n$ steps of the evolution.
\begin{thm} \label{limKnv/Kn}
Let $0<p<1$, $0<q$.
Then for all $w=1,2,\dots$ we have
\begin{equation}
\dfrac{K \left( n,w \right)}{K_n} \rightarrow t_{w}
\end{equation}
almost surely as $n \rightarrow \infty $, where $t_{w}$, $w=1,2, \dots$\,, are positive numbers
satisfying the recurrence relation
\begin{equation} \label{rekurziox(t)-re}
t_{1} = \dfrac{1}{h +1},  \qquad
t_{w} = \dfrac{h \left( w-1 \right) }{h w +1}t_{w-1}, \quad \text{if} \quad w > 1,
\end{equation}
where $h=(1-p)q$.
Moreover, 
\begin{equation} \label{x_{w}_asz2}
t_{w} \sim C  w^{- \left( 1 + \frac{1}{h} \right)}
\end{equation}
as $w \to \infty$,  with $C =\frac{1}{h}\varGamma \left( 1 + \frac{1}{h} \right) $.
\end{thm}
That is the weight distribution of the triangles is asymptotically a power law.

Now, we present the main result of this paper.
Let $E_n$ denote the number of all edges after the $n$th step.
Let $E(n,v)$ denote the number of edges of weight $v$ after the $n$th step.
Let
\begin{equation}  \label{ab}
a= \dfrac{1}{3}pr +  \left(1-p\right) q,  \quad
b = \dfrac{2}{p^2} \left[p \left( 1-r \right) + 3\left( 1-p \right)\left( 1-q \right)\right].
\end{equation}
Next theorem shows the asymptotic power law distribution of the weights of the edges.
%
\begin{thm} \label{limEnv/En}
Let $0<p\le 1$.
Assume that either $r>0$ or $(1-p)q>0$.
Then for any fixed $v$ we have
\begin{equation}   \label{Env/En}
\dfrac{E \left( n,v\right)}{E_n} \rightarrow u_{v}
\end{equation}
almost surely as $n \rightarrow \infty $, where $u_{v}$, $v=1,2,\dots$ are fixed positive numbers.
The numbers $u_{v}$ satisfy the following recurrence relation
\begin{equation} \label{rekurzio_ev}
u_1 = \dfrac{1}{a+1},  \qquad u_v = \dfrac{(v-1)a}{va+1} u_{v-1}, \qquad   v>1.
\end{equation}
The sequence  $u_1, u_2,\dots$ is a proper discrete probability distribution.
Furthermore,
\begin{equation*}
u_v \sim \dfrac{\varGamma \left( 1 + \frac{1}{a} \right)}{a} v^{- \left( 1 + \frac{1}{a} \right)} \quad \text{as} \quad v \to \infty.
\end{equation*}
\end{thm}
The proof of Theorem~\ref{limEnv/En} will be given in Section~\ref{Sect3}.
It is easy to see that Theorem~\ref{limKnv/Kn} is true in a more general setting, that is in the general $N$-interactions model.
Therefore in Section~\ref{Sect4} we introduce the $N$-interactions model.
We describe the asymptotic behaviour of the weights of the $N$-cliques (Theorem~\ref{limKnv/Kn+}).
We can see that Theorem~\ref{limKnv/Kn} is a particular case of Theorem~\ref{limKnv/Kn+} for $N=3$.
We remark that in the general $N$-interactions model the power law distributions both for the weights and degrees of vertices are known  (see \cite{FaPo1} and \cite{FaPo2}).

\section{Proof of Theorem~\ref{limEnv/En} and auxiliary results} \label{Sect3}
\setcounter{equation}{0}
%
\begin{rem}  \label{rem2.1}
Each edge has initial weight $1$. 
As at the initial step we have one triangle, so $E_0=3$, $E(0,1)=3$ and $E(0,v)=0$ if $v>1$.
The weight of an edge can be increased by $0$ or $1$. 
An edge of weight $v$ has taken part in interactions $v$-times.
\end{rem}
Let $\FD_{n-1}$ denote the $\sigma$-algebra of observable events after the $(n-1)$th step.
We compute the conditional expectation of $E(n,v)$ with respect to $\FD_{n-1}$ for $v \geq 1$.

Let
\begin{equation}  \label{p1}
p_1(n,v)= \dfrac{v}{n} a + \dfrac{p^2}{V_{n-1}(V_{n-1}-1)}b.
\end{equation}

The following lemma contains a basic equation of the paper. 
Let $E(i,0)=0$ for any $i=0,1,2,\dots$.
\begin{lem} \label{EE/F}
Let $\delta_{1,v}$ denote the Dirac-delta.
Then for $v\ge 1$ we have
\begin{equation*}
{\E} \{ E(n,v)|\FD_{n-1} \} =
p_1(n, v-1) E(n-1,v-1)  + (1- p_1(n, v)) E(n-1,v) +   
\end{equation*}
\begin{equation} \label{CondExp}
 + \delta_{1,v} \left[ 2p+ p(1-r) \left(1- \dfrac{E_{n-1}}{\binom{V_{n-1}}{2}} \right) +3 (1-p)(1-q)\left(1- \dfrac{E_{n-1}}{\binom{V_{n-1}}{2}} \right) \right].
\end{equation}
\end{lem}
\begin{pf}
The cumulative weight of the edges after $n-1$ steps is $3n$.
The cumulative weight of triangles after $(n-1)$ steps is $n$.
An edge of weight $v$ is included in triangles having cumulative weight $v$.
Moreover, after $(n-1)$ steps, we have the following.
When we choose two vertices uniformly at random, then the probability that a given edge is chosen is 
$
{1}\big/{\binom{V_{n-1}}{2}}
$.
When we choose three vertices randomly, then the probability that a given edge is chosen is 
$$
\frac{V_{n-1}-2}{\binom{V_{n-1}}{3}} = \frac{3}{\binom{V_{n-1}}{2} }.
$$
Therefore the probability that an edge of weight $v$ takes part in interaction at step $n$ is
$$
p \left(r \dfrac{v}{3n} + (1-r) \dfrac{1}{\binom{V_{n-1}}{2}}\right) +
\left(1-p\right) \left(q \dfrac{v}{n} + (1-q) \dfrac{3}{\binom{V_{n-1}}{2}}\right) = p_1(n,v),
$$
where $a$ and $b$ are defined by \eqref{ab}.
These considerations explain terms 
$$p_1(n, v-1) E(n-1,v-1) \quad {\text {and}} \quad (1- p_1(n, v)) E(n-1,v)$$
in \eqref{CondExp}.

To explain the third therm in \eqref{CondExp}, first we remark that it measures the number of new edges, that is edges of weight  $1$.
If a new vertex is born, then two or three new edges are born.
Therefore we have term
$$
2p+ p(1-r) \left(1- \dfrac{E_{n-1}}{\binom{V_{n-1}}{2}} \right).
$$
If no new vertex is born, then new edge can be born when we choose three old vertices uniformly.
So consider the experiment when we choose randomly three vertices of a graph having $y$ vertices and $e$ edges.
Then the expected number of edges connecting the three vertices chosen is 
$$
\dfrac{e(y-2)}{\binom{y}{3}} = \dfrac{3e}{\binom{y}{2}}.
$$
It explains term 
$$
3 (1-p)(1-q)\left(1- \dfrac{E_{n-1}}{\binom{V_{n-1}}{2}} \right)
$$
in \eqref{CondExp}. \hfill $\Box$
\end{pf}
%
\begin{thm} \label{limEnv/np}
Let $0<p\le 1$.
Assume that either $r>0$ or $(1-p)q>0$.
Then for any fixed $v$ we have
\begin{equation}   \label{E/n}
\dfrac{E \left( n,v\right)}{np} \rightarrow e_{v}
\end{equation}
almost surely as $n \rightarrow \infty $, where $e_{v}$, $v=1,2,\dots$, are fixed positive numbers.
Furthermore, the numbers $e_{v}$ satisfy the following recurrence relation
\begin{equation} \label{rekurzio_ev}
e_1 = \dfrac{p(3-r)+ 3 (1-p)(1-q)}{(a+1)p}, \quad
e_v = \dfrac{(v-1)a}{va+1} e_{v-1}, \qquad   v>1.
\end{equation}
\end{thm}
\begin{pf}
Introduce notation
\begin{equation} \label{c_def}
c(n,v) = \prod_{i=1}^{n} \left( 1-p_1(i,v) \right)^{-1}, \quad n \geq 1, \ v \geq 1.
\end{equation}
$c(n,v)$ is an $\FD_{n-1}$-measurable random variable.
We see that $V_n$ is the sum of i.i.d. Bernoulli random variables with success probability $p$.
Therefore, applying the Marcinkiewicz strong law of large numbers to the number of vertices, we have
\begin{equation}  \label{Markinkiewicz}
V_n = pn + \o \left( n^{1/2 + \varepsilon} \right)
\end{equation}
almost surely, for any $\varepsilon >0$.
Equation \eqref{Markinkiewicz} means that
$ \lim_{n\to\infty} (V_n-pn)/(n^{1/2+\varepsilon}) =0$ almost surely. 

Using \eqref{Markinkiewicz} and the Taylor expansion for $\log(1+x)$, we obtain
\begin{equation*}
\log c\left( n,v \right) = 
-\sum_{i=1}^{n} \log  \left( 1-\dfrac{av}{i}-\dfrac{b p^2}{ V_{i-1} \left( V_{i-1} -1 \right)} \right) =
\end{equation*}
\begin{equation*} =
-\sum_{i=1}^{n} \log  \left( 1-\dfrac{av}{i}-\dfrac{b}{i^2 + \o \left( i^{3/2 + \varepsilon} \right)} \right) = 
va  \sum_{i=1}^{n} \dfrac{1}{i} + \O \left( 1 \right),
\end{equation*}
where the error term is convergent as $n \to \infty$.
Therefore
\begin{equation} \label{c(n,v)_asz.}
c(n,v) \sim a_{v} n^{av}
\end{equation}
almost surely as $n \to \infty$, where $a_{v}$ is a positive random variable.

Let
\begin{equation*} \label{Z_def}
Z \left( n,v\right) = c\left( n,v \right) E \left( n,v \right) \quad \text{for} \quad 1 \leq v.
\end{equation*}
Using \eqref{CondExp}, we can see that $ \left\{ Z \left( n,v \right) , \FD_{n} , n=1,2,\dots \right\}$ is a non-negative submartingale for any fixed 
$v \geq 1$.
Applying the Doob-Meyer decomposition to $ Z \left( n,v \right)$, we obtain
\begin{equation*} \label{Doob_Meyer}
Z \left( n,v \right) = M \left( n,v \right) + A \left( n,v \right),
\end{equation*}
where $M \left( n,v \right)$ is a martingale and $A \left( n,v \right)$ is a predictable increasing process. 
The general form of $A \left( n,v \right)$ is the following:
\begin{equation}  \label{A_alt}
A \left( n,v \right) = {\E} Z \left( 1,v \right) + \sum_{i=2}^{n} \left[ {\E}  \left( Z \left( i,v \right) | \FD_{i-1} \right) - Z \left( i-1,v \right)  \right].
\end{equation}
Using \eqref{CondExp} and \eqref{A_alt}, we obtain
\begin{equation} \label{A(n,v)}
A \left( n,v \right)= 
{\E} Z \left( 1,v \right)+
\sum_{i=2}^{n} c \left( i,v \right) \left[ p_1(i,v-1) E \left( i-1,v-1 \right) + \delta_{1,v} Q\right],
\end{equation}
where $Q$ denotes an appropriate term.
Let $B \left( n,v \right)$ be the sum of the conditional variances of the process $Z \left( n,v \right)$.
We give an upper bound of $B \left( n,v \right)$.
\begin{equation*} 
B \left( n,v \right) = \sum_{i=2}^{n} {\D}^2 \left( Z \left( i,v \right) | \FD_{i-1} \right) =
 \sum_{i=2}^{n} {\E} \{ \left( Z \left( i,v \right) - {\E} \left( Z \left( i,v \right) |\FD_{i-1}\right) \right)^2 | \FD_{i-1}  \} = 
\end{equation*}
\begin{equation}\label{Bnv}
=
 \sum_{i=2}^{n} c\left( i,v \right)^2 {\E} \{ \left( E \left( i,v \right) - {\E} \left( E \left( i,v\right) |\FD_{i-1}\right) \right)^2 | \FD_{i-1} 
 \} \leq 
\end{equation}
\begin{equation*}
\leq 
\sum_{i=2}^{n} c\left( i,v \right)^2 {\E} \{ \left( E \left( i,v \right) -  E \left( i-1,v \right) \right)^2 | \FD_{i-1} 
 \} \leq
9 \sum_{i=2}^{n} c\left( i,v \right)^2 = \O\left( n^{2 va +1} \right).
\end{equation*}
Above we used that $c\left( i,v \right)$ is $\FD_{i-1}$-measurable, applied known properties of the conditional expectation, 
and took into account that at each step at most three edges are involved in interaction.

Now we prove \eqref{E/n} by induction.
First let $v=1$. 
By \eqref{A(n,v)}, \eqref{Markinkiewicz}, \eqref{E/n+} and \eqref{c(n,v)_asz.}, we obtain
\begin{align}
A \left( n,1 \right) &={\E} Z \left( 1,1 \right)+\notag\\
&+ \sum_{i=2}^{n} c \left( i,1 \right) 
\left[ 2p+ p(1-r) \left(1- \dfrac{E_{i-1}}{\binom{V_{i-1}}{2}} \right) +3 (1-p)(1-q)\left(1- \dfrac{E_{i-1}}{\binom{V_{i-1}}{2}} \right)\right] \notag \sim\\
& \sim  \sum_{i=2}^{n} a_{1} i^{a} \left[ 2p+ p(1-r) +3 (1-p)(1-q) \right] \notag  \sim\\
&  \sim  a_{1} \dfrac{n^{a+ 1}}{a+ 1} \left[ 2p+ p(1-r) +3 (1-p)(1-q) \right] \label{Teljes indukcio:v=1}
\end{align}
a.s. as $n \to \infty$. 
By \eqref{Bnv}, $B \left( n,1 \right) = \O \left( n^{2a + 1} \right)$ and therefore 
$ \left(B \left( n,1 \right)\right)^{\frac{1}{2}} \log B \left( n,1 \right) = \O \left( A \left( n,1 \right)\right)$.
It follows from Proposition VII-2-4 of \cite{neveu} that
\begin{equation} \label{Z(n,1)_A}
Z \left( n,1 \right) \sim A \left( n,1 \right) \quad \text{a.s. on the event } \quad \{A \left( n,1 \right) \to \infty\} \quad \text{as}\quad n \to \infty.
\end{equation}
As, by \eqref{Teljes indukcio:v=1}, $A(n,1) \to \infty$ a.s., therefore
using \eqref{c(n,v)_asz.}, relation \eqref{Z(n,1)_A} implies
\begin{equation*}
\dfrac{E\left( n,1 \right)}{np} = \dfrac{Z \left( n,1 \right)}{c\left( n,1 \right)np} \sim \dfrac{A \left( n,1 \right)}{c\left( n,1 \right)np} 
\sim \dfrac{ a_{1} \dfrac{n^{a +1}}{a +1}  \left[ 2p+ p(1-r) +3 (1-p)(1-q) \right]  }{a_{1} n^{a }p n} =  
\end{equation*}
\begin{equation*} =
\dfrac{ p(3-r) +3 (1-p)(1-q) }{(a +1)p} = e_{1} > 0
\end{equation*}
almost surely.
So \eqref{E/n} is valid for $v=1$.

Now, let $v>1$ and suppose that the statement \eqref{E/n} is true for $v-1$. 
Then by \eqref{Markinkiewicz}, \eqref{c(n,v)_asz.}, \eqref{A(n,v)} and using the induction hypothesis, we see that
\begin{multline*}
\begin{split}
A \left( n,v \right) &= {\E} Z \left( 1,v \right)+
\sum_{i=2}^{n} c \left( i,v \right)  p_1(i,v-1) E \left( i-1,v-1 \right)  \sim \\
&\sim \sum_{i=2}^{n} c \left( i,v \right)  e_{v-1} ip \left( \frac{a(v-1)}{i} + \frac{b}{i^2}\right)  
 \sim  e_{v-1} p a (v-1) a_v \frac{n^{va+1}}{va+1}.
\end{split}
\end{multline*}
We see that $A(n,v) \to \infty $ a.s.
We also have $\left(B \left( n,v\right)\right)^{\frac{1}{2}} \log B \left( n,v \right) = \O \left( A \left( n,v \right)\right)$.
So, using Proposition VII-2-4 of \cite{neveu}, 
we have $Z \left( n,v \right) \sim A \left( n,v \right)$. 
Therefore
\begin{equation*}
\dfrac{E\left( n,v \right)}{np} = \dfrac{Z \left( n,v \right)}{c\left( n,v \right)np} \sim \dfrac{A \left( n,v \right)}{c\left( n,v \right)np} 
\sim \dfrac{ a_{v}  e_{v-1}  p a (v-1)  \dfrac{n^{av +1}}{va +1} }{a_{v} n^{a v}p n} =  e_{v-1}  \dfrac{a(v-1)}{va +1} =e_v.
\end{equation*}
Therefore \eqref{E/n} is valid.
\hfill$\Box$
\end{pf}
%
The following lemma contains the basic equation for the number of edges. 
\begin{lem} \label{E/F}
\begin{equation} \label{CondExp2}
{\E} \{ E_n|\FD_{n-1} \} = \varrho_{n-1} E_{n-1} +A,
\end{equation}
where
\begin{equation} \label{rhoA}
\varrho_{n-1}= 1- \dfrac{2p(1-r) + 6(1-p)(1-q)}{V_{n-1}(V_{n-1}-1)} \quad {\text{and}} \quad  A= p(3-r) + 3 (1-p) (1-q).
\end{equation}
\end{lem}
\begin{pf}
As the number of edges increases if and only if when the number of edges of weight $1$ increases, therefore by \eqref{CondExp}, we obtain
$$
 {\E} \{ E_n|\FD_{n-1} \} = E_{n-1} +  2p+ p(1-r) \left(1- \dfrac{E_{n-1}}{\binom{V_{n-1}}{2}} \right) +3 (1-p)(1-q)\left(1- \dfrac{E_{n-1}}{\binom{V_{n-1}}{2}} \right).
$$
It implies \eqref{CondExp2}.
\hfill $\Box$
\end{pf}
%
\begin{thm} \label{limEn/n}
Let $0<p\leq 1$.
Then
\begin{equation}   \label{E/n+}
\dfrac{E_n}{n} \rightarrow A
\end{equation}
almost surely as $n \rightarrow \infty $, where $A$ is defined by \eqref{rhoA}.
\end{thm}
\begin{pf}
First we remark that $\varrho_{n}>0$ because $V_n \ge 3$.
Introduce notation $Q_1=1$, 
\begin{equation} \label{c_def}
Q_n = \prod_{i=1}^{n-1} \left( \varrho_{i} \right)^{-1}, \quad n \geq 2.
\end{equation}
Then $Q_n$ is an $\FD_{n-1}$-measurable random variable.
Let
\begin{equation*} \label{Z_def}
Z_n= Q_n E_n \quad \text{for} \quad 1 \leq n.
\end{equation*}
Using \eqref{CondExp2}, we can see that 
\begin{equation} \label{CondExp3}
 {\E} \{ Z_n|\FD_{n-1} \} = Z_{n-1}  + A Q_n.
\end{equation}
Therefore $ \left\{ Z_n, \FD_{n} , n=1,2,\dots \right\}$ is a non-negative submartingale.
Let $A_n$ be the predictable increasing process in the Doob-Meyer decomposition of $ Z_n$.
By \eqref{CondExp3}, we obtain
\begin{equation} \label{A(n)}
A_n= 
{\E} Z_1 + \sum_{i=2}^{n} \left[ {\E} \left( Z_i | \FD_{i-1} \right) - Z_{i-1}  \right]= {\E} Z_1+ A \sum_{i=2}^{n} Q_i.
\end{equation}
Using \eqref{Markinkiewicz} and the Taylor expansion for $\log(1+x)$, we can obtain upper and lower bounds for $Q_i$.
Then
$$
C_1 n <  A_n  < C_2 n,
$$ 
where $C_1$ and $C_2$ are finite positive random variables.
Let $B_n$ be the sum of the conditional variances of the process $Z_n$.
We give an upper bound for $B_n$.
\begin{equation}\label{Bn}
B_n
\leq 
\sum_{i=2}^{n} Q_i^2 {\E} \{ \left( E_i  -  E_{i-1} \right)^2 | \FD_{i-1}  \} \leq
9 \sum_{i=2}^{n} Q_i^2  \leq C_3 n,
\end{equation}
where $C_3$ is a fixed finite random variable.
Above we used that $Q_i$ is $\FD_{i-1}$-measurable and at each step at most three edges are involved in interaction.
Therefore
$ B_n^{\frac{1}{2}} \log B_n  = \O \left( A_n \right)$ almost surely.
It follows from Proposition VII-2-4 of \cite{neveu} that
$
Z_n \sim A_n$ a.s. as $n \to \infty$.
Therefore
\begin{equation*}
\dfrac{E_n}{n} = \dfrac{Z_n}{Q_{n} n} \sim \dfrac{A_n }{Q_{n} n} 
= \dfrac{{\E} Z_1 }{Q_{n} n}   + A \dfrac{1}{ n} \dfrac{1}{Q_{n}}    \sum_{i=2}^{n} Q_{i} \to A  \ \ \ {\text{a.s.}}
\end{equation*}
because the sequence $Q_n$ is positive, increasing and has a finite limit almost surely.
\hfill$\Box$
\end{pf}
{\it {Proof of Theorem \ref{limEnv/En}.}}
By Theorems \ref{limEnv/np}  and \ref{limEn/n},
\begin{equation*}   
\dfrac{E \left( n,v\right)}{E_n}  = \dfrac{p E \left( n,v\right)}{np} \Big/   \dfrac{E_n}{n}   \rightarrow  \dfrac{p e_v}{A} = u_{v}
\end{equation*}
almost surely as $n \rightarrow \infty $.
Furthermore, 
$$
u_1 = \dfrac{ p e_1 }{A} = \dfrac{1}{a+1}.
$$
Multiplying  both sides of \eqref{rekurzio_ev} by $p/A$, we obtain
\begin{equation*} 
u_v = \dfrac{(v-1)a}{va+1} u_{v-1}, \qquad   v>1.
\end{equation*}
Using this equation, we see that
\begin{equation}  \label{uGamma}
u_v = u_1  \dfrac{\varGamma\left( v\right) \varGamma \left( 2 + \frac{1}{a} \right) }{\varGamma \left( 1 + v +\frac{1}{a} \right)}.
\end{equation}
Now, by the Stirling formula,
\begin{equation}  \label{lim_u_v}
u_v \sim u_1  \varGamma \left( 2 + \frac{1}{a} \right)v^{- \left( 1 + \frac{1}{a} \right)} =  
\dfrac{\varGamma \left( 1 + \frac{1}{a} \right)}{a} v^{- \left( 1 + \frac{1}{a} \right)}
\end{equation}
as $v \to \infty$.
We need the following formula (it is included in \cite{prudnikov} and it can be proved by mathematical induction)
$$
\sum_{k=0}^{n} \dfrac{\varGamma \left( k+s \right)}{\varGamma \left( k+t \right)} = 
\dfrac{1}{s-t+1} \left[ \dfrac{\varGamma \left( n+s+1 \right)}{\varGamma \left( n+t \right)} - \dfrac{\varGamma \left( s \right)}{\varGamma \left( t-1 \right)} \right].
$$ 
This formula, \eqref{uGamma} and Stirling's formula imply $ \sum_{v=1}^{n} u_v\to 1$, as $n \to \infty$.
So $\sum_{v=1}^{\infty} u_v =1$. 
Therefore the sequence  $u_1, u_2,\dots$ is a proper discrete probability distribution.
\hfill$\Box$
%
\section{The $N$-cliques in the $N$-interactions model} \label{Sect4}
\setcounter{equation}{0}
Throughout this section we shall study the following $N$-interactions model (see \cite{FaPo2}).
For the sake of brevity, a complete graph with $m$ vertices we call an $m$-clique.
Let $N\ge 3$ be a fixed integer.
At time $n=0$ we start with an $N$-clique.
Let the initial weight of this graph  and the initial weights of its sub-cliques be one. 
(This graph contains $N$ vertices, $ \binom{N}{2}$ edges, \dots , $\binom{N}{m}$ $m$-cliques $\left( m \leq N \right)$. 
Each of these objects has initial weight $1$.)
During the evolution the weight of any new clique is $1$.
At each time step the evolution of the graph is based on the interaction of $N$ vertices. 
The interaction means the following.
At each step $n=1,2,\dots$ we consider $N$ vertices and draw all non-existing edges between those vertices.
So we obtain an $N$-clique.
The weight of this $N$-clique and the weights of all its sub-cliques are increased by $1$. 
(That is we increase the weights of $N$ vertices, $ \binom{N}{2}$ edges, \dots , $N$ different $\left(N-1\right)$-cliques and the $N$-clique itself.)
The selection of the $N$ vertices to interact is the following.

There are two possibilities at each step.
On the one hand, with probability $p$, we add a new vertex and the new vertex and $N-1$ old vertices will interact.
On the other hand, with probability $\left( 1-p \right)$, we do not add any new vertex, but $N$ old vertices interact. 
Here $0 < p \leq 1$ is fixed.

When we add a new vertex, then we choose $N-1$ old vertices and those $N$ vertices will interact.
However, to choose the $N-1$ old vertices we have two options.
With probability $r$ we choose an $\left(N-1\right)$-clique from the existing $\left(N-1\right)$-cliques according to the weights of the $\left(N-1\right)$-cliques.
It means that an $\left(N-1\right)$-clique of weight $w_t$ is chosen with probability $w_t/\sum_h w_h$ (preferential attachment).
On the other hand, with probability $1-r$, we choose $N-1$ out of the existing vertices uniformly, that is all sets consisting of $N-1$ vertices have the same chance.

At a step when we do not add a new vertex, $N$ old vertices interact. 
We have again two options to choose the $N$ old vertices. 
With probability $q$ we choose one $N$-clique out of the existing $N$-cliques according to their weights.
It means that the probability that we choose an $N$-clique is proportional to its weight (preferential attachment).
On the other hand, with probability $1-q$, we choose from the existing vertices uniformly, that is all subsets consisting of $N$ vertices have the same chance.

Let $K(n,w)$ denote the number of $N$-cliques with weight $w$, and let $K_n$ denote the number of all $N$-cliques after $n$ steps of the evolution.
\begin{thm} \label{limKnv/Kn+}
Let $0<p<1$, $0<q$.
Then for all $w=1,2,\dots$ we have
\begin{equation}
\dfrac{K \left( n,w \right)}{K_n} \rightarrow t_{w}
\end{equation}
almost surely as $n \rightarrow \infty $, where $t_{w}$, $w=1,2, \dots$\,, are positive numbers
satisfying the recurrence relation
\begin{equation} \label{rekurziox(t)-re+}
t_{1} = \dfrac{1}{h +1},  \qquad
t_{w} = \dfrac{h \left( w-1 \right) }{h w +1}t_{w-1}, \quad \text{if} \quad w > 1,
\end{equation}
where $h=(1-p)q$.
Moreover, 
\begin{equation} \label{x_{w}_asz2+}
t_{w} \sim C  w^{- \left( 1 + \frac{1}{h} \right)}
\end{equation}
as $w \to \infty$,  with $C =\frac{1}{h}\varGamma \left( 1 + \frac{1}{h} \right) $.
\end{thm}
That is the weight distribution of the $N$-cliques follows asymptotically the same  power law as the weight distribution of the triangles in the three interactions model.
In \cite{FaNo} a heuristic argument was presented to support the above theorem.
\begin{rem}  \label{rem2.1+}
As at the initial step we have one $N$-clique, so $K_0=1$, $K(0,1)=1$ and $K(0,w)=0$ if $w>1$.
At any step the weight of a clique can be increased by $0$ or $1$. 
A clique of weight $w$ has taken part in interactions $w$-times.
\end{rem}
We compute the conditional expectation of $K(n,w)$ with respect to $\FD_{n-1}$ for $w \geq 1$.
Let
\begin{equation}  \label{p1}
p(n,w)= (1-p) \left[q \dfrac{w}{n} + (1-q)\dfrac{1}{\binom{V_{n-1}}{N}} \right].
\end{equation}
\begin{lem} \label{EK/F}
Let $K(n,0)=0$ for any $n$.
For $n,w\ge 1$ we have
\begin{equation*}
{\E} \{ K(n,w)|\FD_{n-1} \} = p(n, w-1) K(n-1,w-1)  + (1- p(n, w)) K(n-1,w) +   
\end{equation*}
\begin{equation} \label{CondExp1}
 + \delta_{1,w} \left[ p+ (1-p)(1-q) \left(1- \dfrac{K_{n-1}}{\binom{V_{n-1}}{N}} \right)\right].
\end{equation}
\end{lem}
\begin{pf}
The probability that the weight of an $N$-clique of weight $w-1$ is increased at step $n$ is $p(n,w-1)$.
The probability that the weight of an $N$-clique of weight $w$ is not increased at step $n$ is $1-p(n,w)$.
The probability that a new $N$-clique is born at step $n$ is 
$$
 p+ (1-p)(1-q) \left(1- \dfrac{K_{n-1}}{\binom{V_{n-1}}{N}} \right).
 $$
 So we obtain \eqref {CondExp1}.
\hfill $\Box$
\end{pf}
\begin{thm} \label{limKnv/n}
Let $0<p<1$ and $0<q$.
For any fixed $w$ we have
\begin{equation}   \label{K/n}
\dfrac{K \left( n,w\right)}{n} \rightarrow k_{w}
\end{equation}
almost surely as $n \rightarrow \infty $, where $k_{w}$, $w=1,2,\dots$, are fixed non-negative numbers.
Furthermore, the numbers $k_{w}$ satisfy the following recurrence relation
\begin{equation} \label{rekurzio_kw}
k_1 = \dfrac{1-h}{h+1}, 
\quad
k_w = \dfrac{(w-1)h}{wh+1} k_{w-1}, \qquad   w>1.
\end{equation}
\end{thm}
\begin{pf}
Introduce notation
\begin{equation} \label{c1_def}
c(n,w) = \prod_{i=1}^{n} \left( 1-p(i,w) \right)^{-1}, \quad n \geq 1, \ w \geq 1.
\end{equation}
$c(n,w)$ is an $\FD_{n-1}$-measurable random variable.
Using \eqref{Markinkiewicz} and the Taylor expansion for $\log(1+x)$, we obtain
\begin{equation*}
\log c\left( n,w \right) = 
\left( wh \right) \sum_{i=1}^{n} \dfrac{1}{i} + \O \left( 1 \right),
\end{equation*}
where the error term is convergent as $n \to \infty$.
Therefore
\begin{equation} \label{c(n,w)_asz.}
c(n,w) \sim h_{w} n^{hw}
\end{equation}
almost surely as $n \to \infty$, where $h_{w}$ is a positive random variable.

Let
\begin{equation*} \label{Z1_def}
Z \left( n,w\right) = c\left( n,w \right) K \left( n,w \right) \quad \text{for} \quad 1 \leq w.
\end{equation*}
Using \eqref {CondExp1}, we can see that $ \left\{ Z \left( n,w \right) , \FD_{n} , n=1,2,,\dots \right\}$ is a non-negative submartingale for any fixed  $w \geq 1$.
We apply the Doob-Meyer decomposition to $ Z \left( n,w \right)$. 
Using \eqref {CondExp1} and \eqref{A_alt}, we obtain that the predictable increasing process is
\begin{equation} \label{A(n,w)}
A \left( n,w \right)= 
{\E} Z \left( 1,w \right)+
\end{equation}
$$ +
\sum_{i=2}^{n} c \left( i,w \right) \left[ p(i,w-1) K \left( i-1,w-1 \right) + \delta_{1,w} \left[ p+ (1-p)(1-q) \left(1- \dfrac{K_{i-1}}{\binom{V_{i-1}}{N}} \right)\right]\right].
$$
Let $B \left( n,w \right)$ be the sum of the conditional variances of the process $Z \left( n,w \right)$.
An upper bound of $B \left( n,w \right)$ is
\begin{equation}\label{Bnw}
B \left( n,w \right) \leq 
 \sum_{i=2}^{n} \left( c\left( i,w \right) \right)^2 = \O\left( n^{2 wh +1} \right).
\end{equation}
Above we used that at each step at most one $N$-clique is involved in interaction.

Now we use induction on $w$.
Let $w=1$. 
By \eqref{A(n,w)} and \eqref{c(n,w)_asz.}, we obtain
\begin{multline} \label{Teljes indukcio:w=1}
\begin{split} 
A \left( n,1 \right) &={\E} Z \left( 1,1 \right)+ 
\sum_{i=2}^{n} c \left( i,1 \right) \left[ p+ (1-p)(1-q) \left(1- \dfrac{K_{i-1}}{\binom{V_{i-1}}{N}} \right)\right] \sim \\
&  \sim  h_{1} \dfrac{n^{h+ 1}}{h+ 1} \left( 1-h \right) 
\end{split}
\end{multline}
a.s. as $n \to \infty$. 
By \eqref{Bnw}, $B \left( n,1 \right) = \O \left( n^{2h + 1} \right)$ and therefore 
$ \left(B \left( n,1 \right)\right)^{\frac{1}{2}} \log B \left( n,1 \right) = \O \left( A \left( n,1 \right)\right)$.
It follows from Proposition VII-2-4 of \cite{neveu} that
\begin{equation} \label{Z(n,1)_A+}
Z \left( n,1 \right) \sim A \left( n,1 \right) \quad \text{a.s. on the event } \quad \{A \left( n,1 \right) \to \infty\} \quad \text{as}\quad n \to \infty.
\end{equation}
As, by \eqref{Teljes indukcio:w=1}, $A(n,1) \to \infty$ a.s., therefore
using \eqref{c(n,w)_asz.}, relation \eqref{Z(n,1)_A+} implies
\begin{equation*}
\dfrac{K\left( n,1 \right)}{n} = \dfrac{Z \left( n,1 \right)}{c\left( n,1 \right)n} \sim \dfrac{A \left( n,1 \right)}{c\left( n,1 \right)n} 
\sim \dfrac{ h_{1} \dfrac{n^{h +1}}{h +1}  \left( 1-h\right)  }{h_{1} n^{h }n} =  
\dfrac{1-h}{h +1} = k_{1} > 0
\end{equation*}
almost surely.
So \eqref{K/n} is valid for $w=1$.

Now, let $w>1$ and suppose that the statement is true for $w-1$. 
Then by \eqref{Markinkiewicz}, \eqref{c(n,w)_asz.}, \eqref{A(n,w)} and using the induction hypothesis, we see that
\begin{multline*}
\begin{split}
A \left( n,w \right) &= {\E} Z \left( 1,w \right)+
\sum_{i=2}^{n} c \left( i,w \right)  p(i,w-1) K \left( i-1,w-1 \right)  \sim \\
&\sim \sum_{i=2}^{n} i^{wh} h_w  k_{w-1} i \left( \frac{h(w-1)}{i} + \frac{(1-p)(1-q)}{\binom{V_{i-1}}{N}}\right)  
 \sim  k_{w-1}  h (w-1) h_w \frac{n^{wh+1}}{wh+1}.
\end{split}
\end{multline*}
We see that $A(n,w) \to \infty $ a.s.
We also have $\left(B \left( n,w\right)\right)^{\frac{1}{2}} \log B \left( n,w \right) = \O \left( A \left( n,w \right)\right)$.
So, using Proposition VII-2-4 of \cite{neveu}, 
we have $Z \left( n,w \right) \sim A \left( n,w \right)$. 
Therefore
\begin{equation*}
\dfrac{K\left( n,w \right)}{n} = \dfrac{Z \left( n,w \right)}{c\left( n,w \right)n} \sim \dfrac{A \left( n,w \right)}{c\left( n,w \right)n} 
\sim \dfrac{ k_{w-1}  h (w-1) h_w \frac{n^{wh+1}}{wh+1} }{h_{w} n^{h w} n} =  k_{w-1}  \dfrac{h(w-1)}{wh +1} =k_w.
\end{equation*}
Therefore \eqref{K/n} is valid.
\hfill$\Box$
\end{pf}
The following lemma contains the basic equation for the number of $N$-cliques. 
\begin{lem} \label{K/F}
\begin{equation} \label{CondExp2+}
\E\{K_n | \FD_{n-1} \} = \varrho_{n-1} K_{n-1} +B,
\end{equation}
where
\begin{equation} \label{rhoB}
\varrho_{n-1}=  1- (1-p)(1-q) \left(1- \dfrac{1}{\binom{V_{n-1}}{N}} \right)\quad {\text{and}} \quad  B= 1-h.
\end{equation}
\end{lem}
\begin{pf}
As the number of $N$-cliques increases if and only if when the number of $N$-cliques of weight $1$ increases, therefore
\eqref{CondExp1} implies
$$
\E\{K_n | \FD_{n-1} \} = K_{n-1} +p + (1-p)(1-q) \left(1- \dfrac{K_{n-1}}{\binom{V_{n-1}}{N}} \right).
$$
So we obtain \eqref{CondExp2+}.
\hfill $\Box$
\end{pf}
%
\begin{thm} \label{limKn/n}
Let $0<p\leq 1$.
Then
\begin{equation}   \label{K/n+}
\dfrac{K_n}{n} \rightarrow B
\end{equation}
almost surely as $n \rightarrow \infty $, where $B$ is defined by \eqref{rhoB}.
\end{thm}
\begin{pf}
First we remark that $\varrho_{n}>0$ because $V_n \ge N$.
Introduce notation
\begin{equation} \label{c1_def}
Q_n = \prod_{i=1}^{n-1} \left( \varrho_{i} \right)^{-1}, \quad n \geq 1.
\end{equation}
Then $Q_n$ is an $\FD_{n-1}$-measurable random variable.
Let
\begin{equation*} \label{Z1_def}
Z_n= Q_n K_n \quad \text{for} \quad 1 \leq n.
\end{equation*}
Using \eqref{CondExp2+}, we see that 
\begin{equation} \label{CondExp3+}
 {\E} \{ Z_n|\FD_{n-1} \} = Z_{n-1}  + B Q_n.
\end{equation}
Therefore $ \left\{ Z_n, \FD_{n} , n=1,2,,\dots \right\}$ is a non-negative submartingale.
Let  $A_n$ be the predictable increasing process in the Doob-Meyer decomposition of $ Z_n$.
Then, by \eqref{A_alt} and \eqref{CondExp3+}, we obtain
\begin{equation} \label{A(n)+}
A_n= {\E} Z_1+ B \sum_{i=2}^{n} Q_i.
\end{equation}
Using \eqref{Markinkiewicz} and the Taylor expansion for $\log(1+x)$, we can obtain upper and lower bounds for $Q_i$.
Then
$$
C_1 n <  A_n  < C_2 n,
$$ 
where $C_1$ and $C_2$ are finite positive constants.
Let $B_n$ be the sum of the conditional variances of the process $Z_n$.
We have the following upper bound for $B_n$
\begin{equation}\label{Bn+}
B_n \leq 
\sum_{i=2}^{n} Q_i^2  \leq C_3 n,
\end{equation}
where $C_3$ is a fixed finite constant.
Above we used that $Q_i$ is $\FD_{i-1}$-measurable and at each step at most one $N$-clique is involved in interaction.
Therefore
$ B_n^{\frac{1}{2}} \log B_n  = \O \left( A_n \right)$.
It follows from Proposition VII-2-4 of \cite{neveu} that
\begin{equation} \label{Z(n)_A}
Z_n \sim A_n \quad \text{a.s. on the event } \quad \{A_n \to \infty\} \quad \text{as}\quad n \to \infty.
\end{equation}
As $A_n \to \infty$ a.s., therefore
\begin{equation*}
\dfrac{K_n}{n} = \dfrac{Z_n}{Q_{n} n} \sim \dfrac{A_n }{Q_{n} n} 
= \dfrac{{\E} Z_1 }{Q_{n} n}   + B \dfrac{1}{ n} \dfrac{1}{Q_{n}}    \sum_{i=2}^{n} Q_{i} \to B  \ \ \ {\text{a.s.}}
\end{equation*}
because the sequence $Q_n$ is positive, increasing and has a finite limit.
\hfill$\Box$
\end{pf}
{\it {Proof of Theorem \ref{limKnv/Kn+}.}}
By Theorems \ref{limKnv/n}  and \ref{limKn/n},
\begin{equation*}   
\dfrac{K \left( n,w\right)}{K_n}  = \dfrac{ K \left( n,w\right)}{n} \Big/   \dfrac{K_n}{n}   \rightarrow  \dfrac{k_w}{B} = t_{w}
\end{equation*}
almost surely as $n \rightarrow \infty $.
Furthermore, 
$$
t_1 = \dfrac{ k_1 }{B} = \dfrac{1}{h+1}.
$$
Dividing  both sides of \eqref{rekurzio_kw} by $B$, we obtain
\begin{equation*} 
t_w = \dfrac{(w-1)h}{wh+1} t_{w-1}, \qquad   w>1.
\end{equation*}
Using this equation, we see that
\begin{equation*}
t_w = \dfrac{1}{h} \dfrac{\varGamma\left( w\right) \varGamma \left( 1 + \frac{1}{h} \right) }{\varGamma \left( 1 + w +\frac{1}{h} \right)}.
\end{equation*}
Now, by the Stirling formula,
\begin{equation}  \label{lim_t_w}
t_w \sim 
\dfrac{1}{h} \varGamma \left( 1 + \frac{1}{h} \right) w^{- \left( 1 + \frac{1}{h} \right)}
\end{equation}
as $w \to \infty$.
We can see that $\sum_{w=1}^{\infty} t_w =1$. 
Therefore the sequence  $t_1, t_2,\dots$ is a proper discrete probability distribution.
\hfill$\Box$
\section{Numerical results} \label{Sect5}
\setcounter{equation}{0}
Here we present some numerical results.
The $3$-interactions model was generated with parameters $p=0.5$, $q=0.5$,  $r=0.5$.
The number of evolution steps was $n=10^7$.
We applied log-log scale to visualize scale-free property.
The diamonds show the simulation results.
The solid lines come from the recursive relations.
The dotted lines show the slope (that is the exponent of the power law).
We obtained that the simulation results fit well to the solid lines which support our theoretical results.
The degree distribution of the vertices is given on Figure~\ref{fokszam}, the weight distribution of the vertices on Figure~\ref{csucssuly},
the weight distribution of the edges on Figure~\ref{elsuly} and finally the weight distribution of the triangles is presented on Figure~\ref{klikksuly}.
Figures \ref{csucssuly}, \ref{elsuly} and \ref{klikksuly} offer numerical evidence for Theorems \ref{theorem:scalefreeWeights}, \ref{limKnv/Kn} and \ref{limEnv/En}, respectively.
Figure~\ref{fokszam} corresponds to Theorem 4.3 of \cite{BaMo2} where the power law degree distribution was presented.
\begin{figure}[!h]
\centering
\includegraphics[width=8cm]{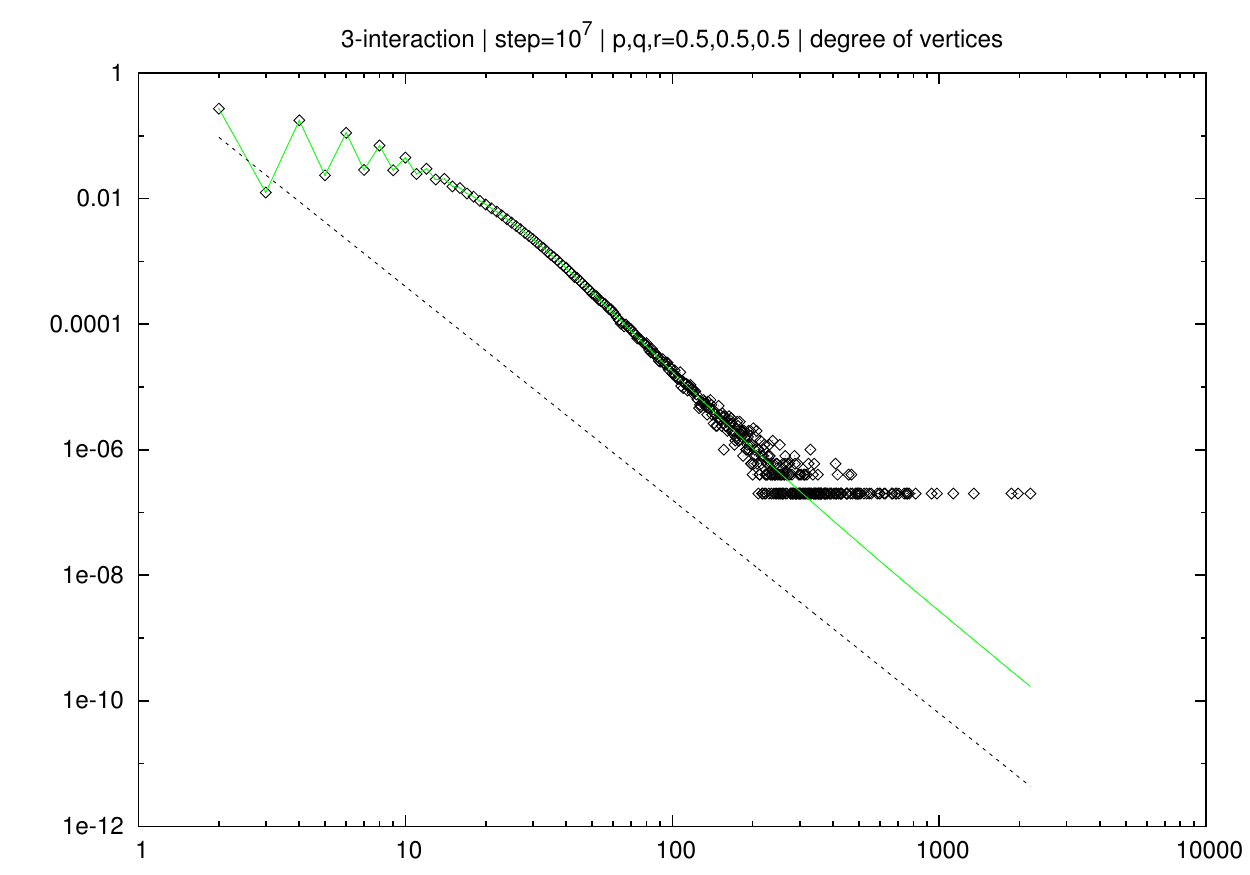}
\caption{The degree distribution of the vertices}
\label{fokszam}
\end{figure}
\begin{figure}[!h]
\centering
\includegraphics[width=8cm]{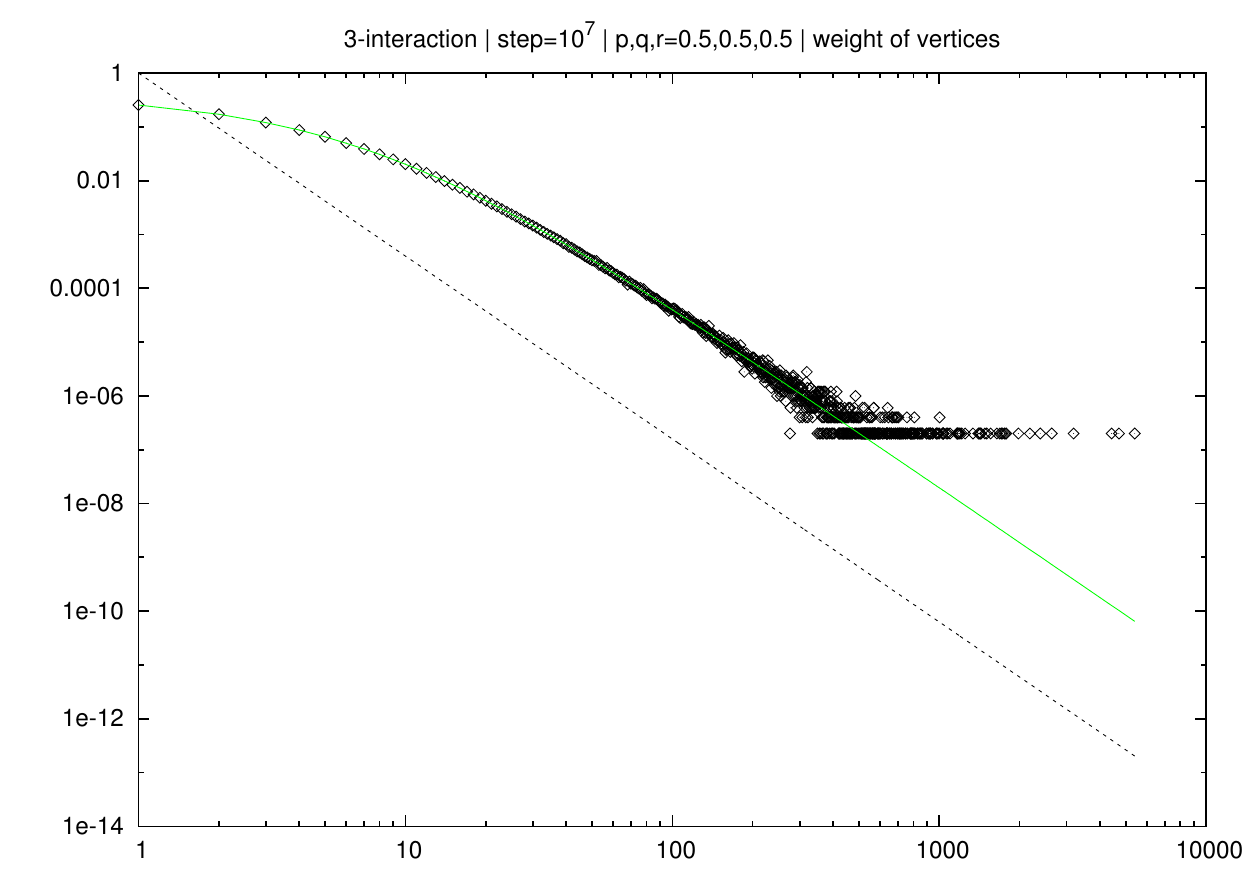}
\caption{The weight distribution of the vertices}
\label{csucssuly}
\end{figure}
\begin{figure}[!h]
\centering
\includegraphics[width=8cm]{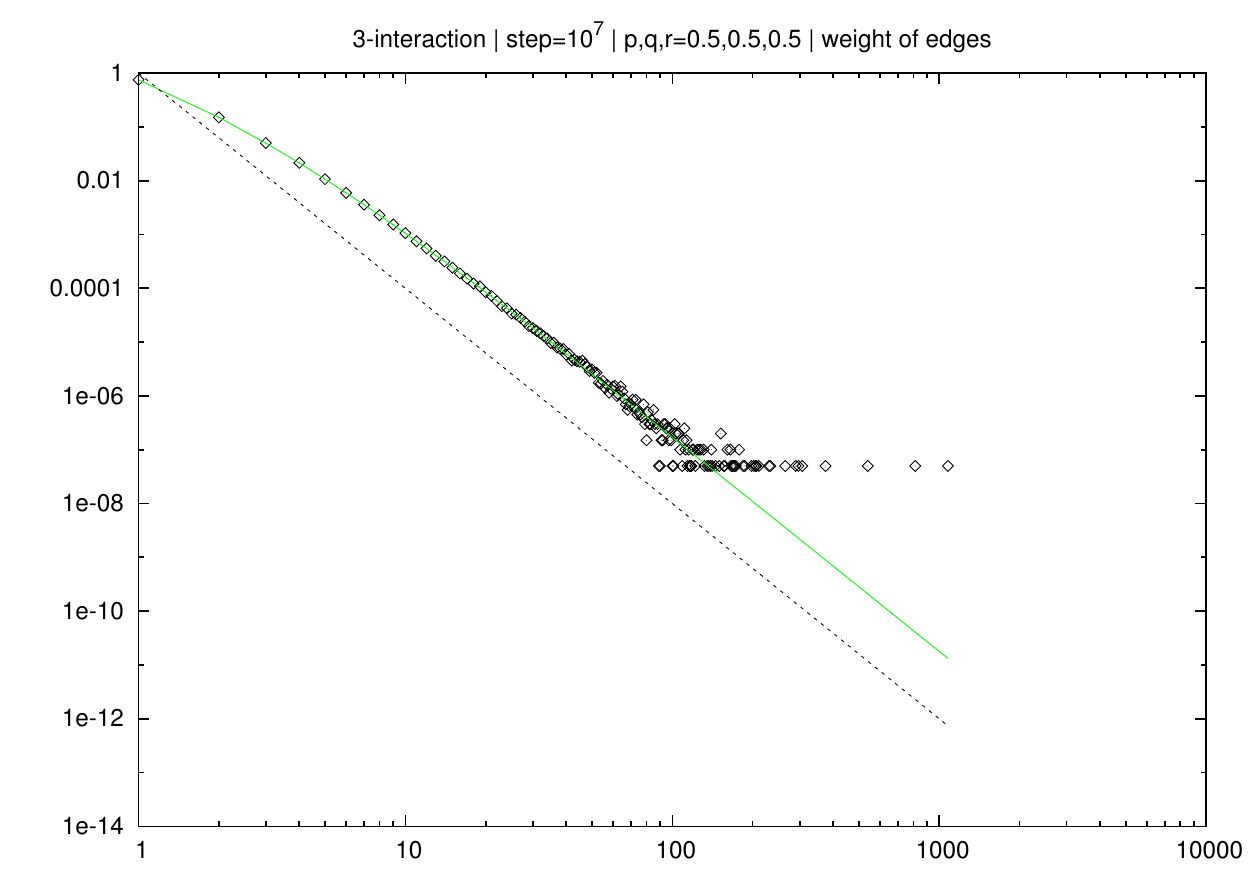}
\caption{The weight distribution of the edges}
\label{elsuly}
\end{figure}
\begin{figure}[!h]
\centering
\includegraphics[width=8cm]{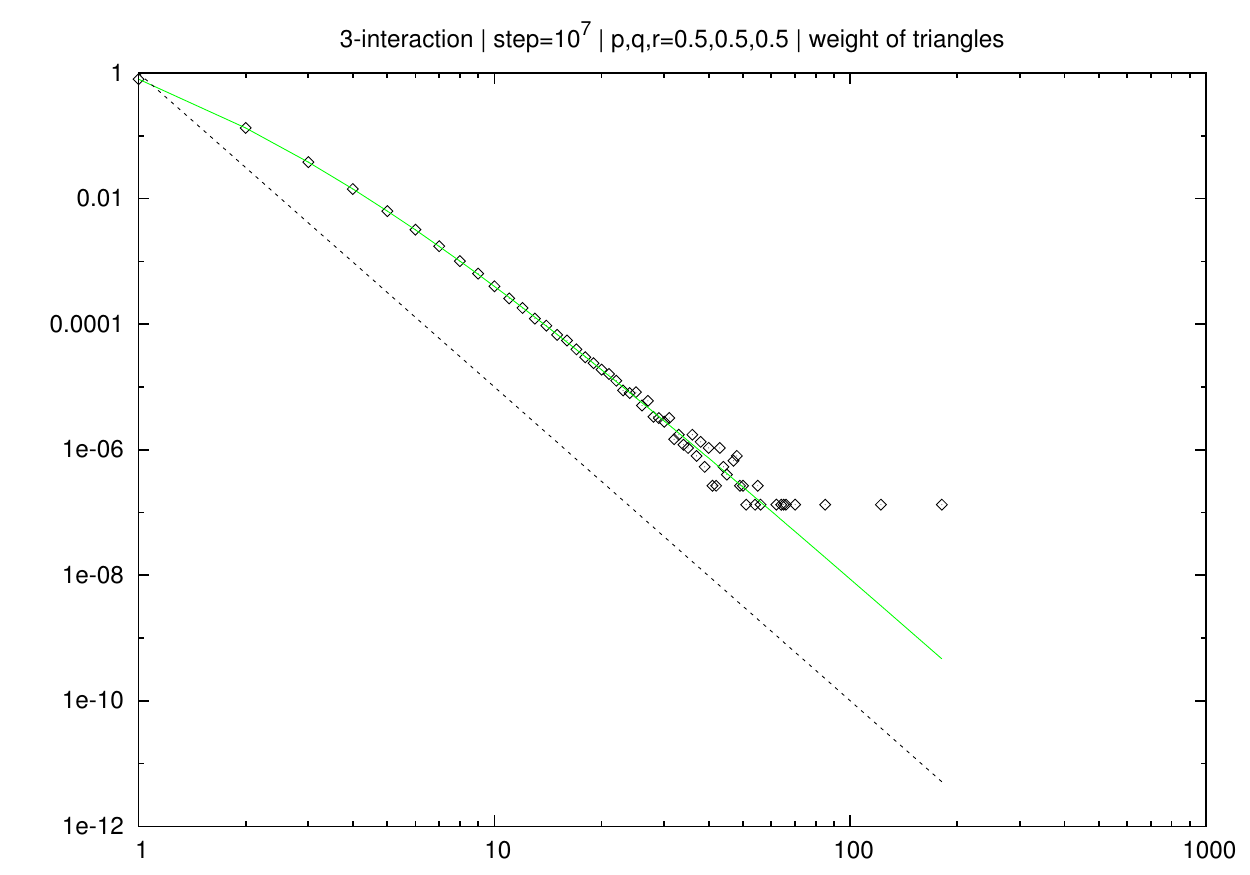}
\caption{The weight distribution of the triangles}
\label{klikksuly}
\end{figure}
%

{\bf Acknowledgement.} We thank the referee for contributing many comments that helped us to correct some mistakes and to clarify the exposition.

\end{document}